\documentclass[12pt,reqno]{amsart}
\usepackage[T1]{fontenc}
\usepackage[utf8]{inputenc}
\usepackage{amssymb, amsthm, amsxtra, epsfig, amsmath, parskip}
\usepackage{hyperref}

\newcommand{\N}{\mathbb{N}}

\newcommand{\E}{\mathbb{E}}

\newcommand{\Prob}{\mathbb{P}}

\newtheorem{theorem}{Theorem}[section]
\newtheorem{lemma}[theorem]{Lemma}

\theoremstyle{definition}

\theoremstyle{remark}

\begin{document}

\sloppy
\title[Upper bound of the constant for the near best m-term approximation]{Derivation of an upper bound of the constant in the error bound for a near best m-term approximation}

\author{Wolfgang Karcher}
\address{Wolfgang Karcher, Ulm University, Institute of Stochastics, Helmholtzstr. 18, 89081 Ulm, Germany}
\email{wolfgang.karcher\@@{}uni-ulm.de}
\author{Hans-Peter Scheffler}
\address{Hans-Peter Scheffler, University of Siegen, Fachbereich 6, Mathematik, Emmy-Noether-Campus, Walter-Flex-Str. 3, 57068 Siegen, Germany}
\email{scheffler\@@{}mathematik.uni-siegen.de}
\author{Evgeny Spodarev}
\address{Evgeny Spodarev, Ulm University, Institute of Stochastics, Helmholtzstr. 18, 89081 Ulm, Germany}
\email{evgeny.spodarev\@@{}uni-ulm.de}

\date{7 October 2009}

\begin{abstract}
In \cite{Tem98}, Temlyakov provides an error bound for a near best $m$-term approximation of a function $g \in L^p([0,1]^d)$, $1<p<\infty$, $d \in \N$, using a basis $L^p$-equivalent to the Haar system $\mathcal{H}$. The bound includes a constant $C(p)$ that is not given explicitly. The goal of this paper is to find an upper bound of the constant for the Haar system $\mathcal{H}$, following the proof in \cite{Tem98}.
\end{abstract}


\maketitle

\baselineskip=18pt

\section{Determining the constant in the one-dimensional case}

Let $\mathcal{H} := \{H_I\}_I$ be the Haar basis in $L^p[0,1]$ indexed by dyadic intervals \\ $I=[(j-1)2^{-n},j 2^{-n})$, $j=1,...,2^n$, $n=0,1,...$ and $I = [0,1]$ with
\begin{eqnarray*}
 H_{[0,1]}(x) &=& 1 \quad \text{for } x \in [0,1), \\
 H_{[(j-1)2^{-n},j2^{-n})}(x) &=& \begin{cases}
                                 2^{n/2}, & x \in [(j-1)2^{-n}, (j-\frac{1}{2})2^{-n}), \\
				 -2^{-n/2}, & x \in [(j-\frac{1}{2})2^{-n}, j2^{-n}), \\
				 0, & \text{otherwise.}
                                \end{cases}
\end{eqnarray*}

Let
$$f = \sum_I c_I(f) H_I,$$
where
$$ c_I(f) := (f,H_I) = \int_0^1 f(x) H_I(x) dx,$$
and denote
$$c_I(f,p) := \Vert c_I(f) H_I \Vert_p.$$

Then $c_I(f,p) \to 0$ as $|I| \to 0$.

Denote by $\Lambda_m$ a set of $m$ dyadic intervals $I$ such that
$$\min\limits_{I \in \Lambda_m} c_I(f,p) \geq \max\limits_{J \notin \Lambda_m} c_J(f,p).$$
This means that $\Lambda_m$ contains the $m$ largest values of $c_I(f,p)$ where $I$ runs through all dyadic intervals. Then we define the Greedy algorithm $G_m^p(\cdot, \mathcal{H})$ as
$$G_m^p(f, \mathcal{H}) := \sum_{I \in \Lambda_m} c_I(f) H_I.$$

The following theorem provides an error bound for the approximation of a function $f \in L^p[0,1]$ by the Greedy algorithm $G_m^p(\cdot, \mathcal{H})$:

\begin{theorem} \label{th:2}
 Let $1 < p < \infty$. Then for any $g \in L^p[0,1]$, we have
 $$\Vert g - G_m^{p} (g, \mathcal{H}) \Vert_p \leq \left(2+\frac{1}{\left(1-\left(\frac{1}{2}\right)^{1/p}\right)^2}\right) \cdot \left( \max\left(p,\frac{p}{p-1}\right) -1\right)^2 \cdot \sigma_m(g)_p. $$
\end{theorem}

\begin{proof}
 The Littlewood-Paley theorem for the Haar system gives for $1<p<\infty$
\begin{equation}
 C_3(p) \left\Vert \left( \sum_I |c_I(g) H_I|^2 \right)^{\frac{1}{2}} \right\Vert_p \leq ||g||_p \leq C_4(p) \left\Vert \left( \sum_I |c_I(g) H_I|^2 \right)^{\frac{1}{2}} \right\Vert_p. \label{eq:littlewood}
\end{equation}

In case of $g$ being a martingale, explicit formulas for these constants are known (cf.~\cite{Bur88}). In Lemma \ref{lemma:cond_mart1}, page \pageref{lemma:cond_mart1}, it is shown that the Haar series
$$ g = \sum_I c_I(g) H_I $$
is in fact a (conditionally symmetric) martingale. 

Thus, taking the constants in \cite{Bur88}, page 87, we have
$$C_3(p) = \frac{1}{\max\left(p,\frac{p}{p-1}\right)-1} \quad \text{and} \quad C_4(p) = \max\left(p,\frac{p}{p-1}\right)-1. $$

Let $T_m$ be an $m$-term Haar polynomial of best $m$-term approximation to $g$ in $L^p[0,1]$:
$$T_m = \sum_{I \in \Lambda} a_I H_I, \quad |\Lambda| = m.$$
For any finite set $Q$ of dyadic intervals we denote by $S_Q$ the projector
$$S_Q(f) := \sum_{I \in Q} c_I(f) H_I.$$
With these definitions, one can derive the following inequality:
\begin{eqnarray*}
 \Vert g - S_\Lambda(g) \Vert_p &=& \Vert g - T_m - S_{\Lambda}(g-T_m) \Vert_p \\
				&\leq& \Vert Id - S_\Lambda \Vert_{p \to p} \sigma_m(g)_p \\
				&\leq& C_4(p) C_3(p)^{-1} \sigma_m(g)_p,
\end{eqnarray*}
where $Id$ denotes the identical operator. The last inequality holds since $$\Vert g \Vert_p \leq 1$$ implies by the Littlewood-Paley theorem (cf.~$(\ref{eq:littlewood})$) that $$\left\Vert \left( \sum_I |c_I(g) H_I|^2 \right)^{\frac{1}{2}} \right\Vert_p \leq C_3(p)^{-1},$$ so that by again applying $(\ref{eq:littlewood})$ we get
\begin{eqnarray*}
 \Vert (Id - S_\Lambda)(g) \Vert_p &=&\left\Vert \sum\limits_{I \notin \Lambda} c_I(g) H_I \right\Vert_p \\
	&\leq& C_4(p) \left\Vert \left( \sum_{I \notin \Lambda} |c_I(g) H_I|^2 \right)^{\frac{1}{2}} \right\Vert_p \\
	&\leq& C_4(p) \left\Vert \left( \sum_{I} |c_I(g) H_I|^2 \right)^{\frac{1}{2}} \right\Vert_p \\
	&\leq& C_4(p) C_3(p)^{-1}.
\end{eqnarray*}

With $C_3(p)$ and $C_4(p)$ given above we have
\begin{equation}
 \Vert g - S_\Lambda(g) \Vert_p \leq \left( \max\left(p,\frac{p}{p-1}\right) -1\right)^2 \cdot \sigma_m(g)_p. \label{eq:temp1}
\end{equation}

Since
$$G_m^p(g) = S_{\Lambda_m}(g),$$
we have
\begin{eqnarray}
 &&\Vert g - G_m^p(g) \Vert_p \nonumber\\
	&\leq& \Vert g - S_\Lambda(g) \Vert_p + \Vert S_\Lambda(g) - S_{\Lambda_m}(g) \Vert_p\nonumber \\
	&\overset{\textup{(\ref{eq:temp1})}}{\leq}& \left( \max\left(p,\frac{p}{p-1}\right) -1\right)^2 \cdot \sigma_m(g)_p + \Vert S_\Lambda(g) - S_{\Lambda_m}(g) \Vert_p. \label{eq:temp2}
\end{eqnarray}

It remains to estimate $\Vert S_\Lambda(g) - S_{\Lambda_m}(g) \Vert_p$ appropriately:
\begin{eqnarray}
 \Vert S_\Lambda(g) - S_{\Lambda_m}(g)\Vert_p &=& \Vert S_{\Lambda \setminus \Lambda_m}(g) - S_{\Lambda_m \setminus \Lambda}(g)\Vert_p \nonumber \\
 &\leq& \Vert S_{\Lambda \setminus \Lambda_m}(g) \Vert_p + \Vert S_{\Lambda_m \setminus \Lambda}(g)\Vert_p. \label{eq:temp3}
\end{eqnarray}

The second term in the last expression can be estimated by
\begin{eqnarray}
 \Vert S_{\Lambda_m \setminus \Lambda}(g) \Vert_p &=& \Vert \left(Id - S_{(\Lambda_m \setminus \Lambda)^C}\right) (g) \Vert_p \nonumber\\
	&=& \Vert g - T_m - S_{\Lambda \bigcup \Lambda_m^C}(g-T_m) \Vert_p \nonumber\\
	&\leq& \Vert Id - S_{\Lambda \bigcup \Lambda_m^C} \Vert_{p \to p} \sigma_m(g)_p\nonumber \\
	&\leq& \frac{C_4(p)}{C_3(p)} \sigma_m(g)_p\nonumber \\
	&=& \left( \max\left(p,\frac{p}{p-1}\right) -1\right)^2 \cdot \sigma_m(g)_p. \label{eq:temp4}
\end{eqnarray}

Furthermore
\begin{equation}
 \Vert S_{\Lambda \setminus \Lambda_m}(g) \Vert_p \leq \frac{1}{\left(1-\left(\frac{1}{2}\right)^{1/p}\right)^2} \cdot \Vert S_{\Lambda_m \setminus \Lambda}(g) \Vert_p \label{eq:temp5}
\end{equation}

which will be derived in the following lemmas (Lemma \ref{lemma:1} - \ref{lemma:4}).

Combining (\ref{eq:temp2})--(\ref{eq:temp5}), we get
$$\Vert g - G_m^p(g) \Vert_p \leq \left(2+\frac{1}{\left(1-\left(\frac{1}{2}\right)^{1/p}\right)^2}\right) \cdot \left( \max\left(p,\frac{p}{p-1}\right) -1\right)^2 \cdot \sigma_m(g)_p. $$
\end{proof}

 \begin{lemma} \label{lemma:1}
  Let $n_1 < n_2 < \cdots < n_s$ be integers and let $E_j \subset [0,1]$ be measurable sets, $j=1,...,s$. Then for any $0 < q < \infty$ we have
 \begin{equation}
  \int_0^1 \left(\sum_{j=1}^s 2^{n_j/q} \chi_{E_j}(x) \right)^q dx \leq \left(\frac{1}{1-\left(\frac{1}{2}\right)^{1/q}}\right)^q \cdot \sum_{j=1}^s 2^{n_j} |E_j|. \nonumber
 \end{equation}
 where $\chi_I(\cdot)$ is the characteristic function of the interval $I$:
 $$\chi_I(x) = \begin{cases}
                1, & x \in I, \\
		0, & x \notin I.
               \end{cases}$$
 \end{lemma}
 \begin{proof}
  Denote
	$$F(x) := \sum_{j=1}^s 2^{n_j/q} \chi_{E_j}(x)$$
  and estimate it on the sets
  $$E_l^- := E_l \setminus \bigcup_{k=l+1}^s E_k, \quad l=1,...,s-1; \quad E_s^- := E_s.$$
  We have for $x \in E_l^-$
  \begin{eqnarray*}
   F(x) &\leq& \sum_{j=1}^l 2^{n_j/q} \\
	&=& 2^{n_l/q} \left( \frac{2^{n_1/q}}{2^{n_l/q}} + \cdots +1\right) \\
	&\leq& 2^{n_l/q} \sum_{i=0}^\infty \left(\frac{1}{2^{1/q}}\right)^i \\
	&=& 2^{n_l/q} \frac{1}{1-\left(\frac{1}{2}\right)^{1/q}}.
  \end{eqnarray*}
  Therefore,
 $$ \int_0^1 F(x)^q dx \leq \left(\frac{1}{1-\left(\frac{1}{2}\right)^{1/q}} \right)^q \sum_{l=1}^s 2^{n_l} |E_l^-| \leq \left(\frac{1}{1-\left(\frac{1}{2}\right)^{1/q}} \right)^q \sum_{l=1}^s 2^{n_l} |E_l|,$$
which proves the lemma.
\end{proof}
 
 \begin{lemma} \label{lemma:2}
  Consider
	$$f = \sum_{I \in Q} c_I H_I, \quad |Q| = N.$$
 Let $1 \leq p < \infty$. Assume that
	\begin{equation}
	 \Vert c_I H_I \Vert_p \leq 1, \quad I \in Q. \label{eq:assumption1}
	\end{equation}
 Then
	$$\Vert f \Vert_p \leq \frac{1}{1-\left(\frac{1}{2}\right)^{1/p}} N^{1/p}.$$
 \end{lemma}

 \begin{proof}
  Denote by $n_1 < n_2 < \cdots < n_s$ all integers such that there is $I \in Q$ with $|I| = 2^{-n_j}$. Introduce the sets
  $$E_j := \bigcup_{I \in Q: |I| = 2^{-n_j}} I.$$
  Then the number $N$ of elements in $Q$ can be written in the form
  $$N = \sum_{j=1}^s |E_j| 2^{n_j}.$$
  Furthermore, we have
	$$\Vert c_I H_I\Vert_p = |c_I| |I|^{1/p - 1/2}.$$
  The assumption (\ref{eq:assumption1}) implies $ |c_I| \leq |I|^{1/2 - 1/p}$.
  Next, we have
	$$\Vert f \Vert_p \leq \left\Vert \sum_{I \in Q} |c_I H_I| \right\Vert_p \leq \left\Vert \sum_{I \in Q} |I|^{-1/p} \chi_I(x) \right\Vert_p.$$
  The right hand side of this inequality cna be rewritten as
  $$Y := \left( \int_0^1 \left( \sum_{j=1}^s 2^{n_j/p} \chi_{E_j}(x) \right)^p dx \right)^{1/p}.$$
  Applying Lemma \ref{lemma:1} with $q=p$, we get
  $$\Vert f \Vert_p \leq Y \leq \frac{1}{1-\left(1/2\right)^{1/p}} \left(\sum_{j=1}^s |E_j| 2^{n_j} \right)^{1/p} = \frac{1}{1-\left(1/2\right)^{1/p}} N^{1/p}.$$
\end{proof}

 \begin{lemma} \label{lemma:3}
  Consider
	$$f = \sum_{I \in Q} c_I H_I, \quad |Q| = N.$$
 Let $1 \leq p < \infty$. Assume
	\begin{equation*}
	 \Vert c_I H_I \Vert_p \geq 1, \quad I \in Q.
	\end{equation*}
 Then
	$$\Vert f \Vert_p \geq \left(1-\left(\frac{1}{2}\right)^{1/p}\right) N^{1/p}.$$
 \end{lemma}
 \begin{proof}
  Define
  $$ u := \sum_{I \in Q} \overline{c}_I |c_I|^{-1} |I|^{1/p-1/2} H_I,$$
  where the bar means complex conjugate number. Then for $p' = \frac{p}{p-1}$ we have
  $$\Vert \overline{c}_I |c_I|^{-1} |I|^{1/p-1/2} H_I \Vert_{p'} = 1$$
  and, by Lemma \ref{lemma:2}
  $$\Vert u \Vert_{p'} \leq \frac{1}{1-\left(\frac{1}{2}\right)^{1/p}} N^{1/p'}.$$
  Consider $(f,u)$. We have on the one hand
  $$(f,u)  = \sum_{I \in Q} |c_I| |I|^{1/p-1/2} = \sum_{I \in Q} \Vert c_I H_I \Vert_p \geq N,$$
  and on the other hand
  $$ (f,u)  \leq \Vert f \Vert_p \Vert u \Vert_{p'},$$
  so that
  $$N \leq (f,u)  \leq \Vert f \Vert_p \Vert u \Vert_{p'} \leq \Vert f \Vert_p \frac{1}{1-\left(\frac{1}{2}\right)^{1/p}} N^{1/p'}$$
  which implies
  $$\Vert f \Vert_p \geq \left(1-\left(\frac{1}{2}\right)^{1/p}\right) N^{1/p}.$$
\end{proof}

\begin{lemma} \label{lemma:4}
 Let $1 < p < \infty$. Then for any $g \in L^p[0,1]$ we have
 $$\Vert S_{\Lambda \setminus \Lambda_m}(g) \Vert_p \leq \frac{1}{\left(1-\left(\frac{1}{2}\right)^{1/p}\right)^2} \cdot \Vert S_{\Lambda_m \setminus \Lambda}(g) \Vert_p.$$
\end{lemma}

\begin{proof}
 Denote
 $$A := \max\limits_{I \in \Lambda \setminus \Lambda_m} \Vert c_I(g) H_I \Vert_p \quad \text{and} \quad B:= \min\limits_{I \in \Lambda_m \setminus \Lambda} \Vert c_I(g) H_I \Vert_p.$$
 Then by the definition of $\Lambda_m$ we have
 $$ B \geq A.$$
 Using Lemma \ref{lemma:2}, we get
 \begin{equation}
  \Vert S_{\Lambda \setminus \Lambda_m}(g)\Vert_p \leq A \cdot \frac{1}{1-\left(\frac{1}{2}\right)^{1/p}} \cdot|\Lambda \setminus \Lambda_m |^{1/p} \leq B \cdot \frac{1}{1-\left(\frac{1}{2}\right)^{1/p}} \cdot|\Lambda \setminus \Lambda_m|^{1/p}. \label{eq:temp6}
 \end{equation}

 Using Lemma \ref{lemma:3}, we get
 $$\Vert S_{\Lambda_m \setminus \Lambda}(g)\Vert_p \geq B \cdot \left(1-\left(\frac{1}{2}\right)^{1/p}\right) \cdot |\Lambda_m \setminus \Lambda|^{1/p}$$
 so that
 \begin{equation}
  |\Lambda_m \setminus \Lambda|^{1/p} \leq \frac{1}{B \cdot \left(1-\left(\frac{1}{2}\right)^{1/p}\right)} \Vert S_{\Lambda_m \setminus \Lambda}(g)\Vert_p. \label{eq:temp7}
 \end{equation}
 Since $\vert \Lambda \vert = \vert \Lambda_m \vert = m$, we have $|\Lambda_m \setminus \Lambda| = |\Lambda \setminus \Lambda_m|$ and finally get
 \begin{eqnarray*}
  \Vert S_{\Lambda \setminus \Lambda_m}(g)\Vert_p &\overset{\textup{(\ref{eq:temp6})}}{\leq}& B \cdot \frac{1}{1-\left(\frac{1}{2}\right)^{1/p}} \cdot|\Lambda \setminus \Lambda_m|^{1/p} \\
	&\overset{\textup{(\ref{eq:temp7})}}{\leq}& \frac{1}{\left(1-\left(\frac{1}{2}\right)^{1/p}\right)^2} \Vert S_{\Lambda_m \setminus \Lambda}(g)\Vert_p.
 \end{eqnarray*}
\end{proof}

\begin{lemma} \label{lemma:cond_mart1}
 Let $f \in L^p[0,1]$. Then the Haar series
 $$ g = \sum_I c_I(g) H_I $$
 is a conditionally symmetric martingale.
\end{lemma}

\begin{proof}
First, we give a definition of a conditionally symmetric martingale (cf.~\cite{Wan91a} and \cite{Wan91b}).

Let $(\Omega, \mathcal{F}, \Prob)$ be a probability space with a nondecreasing sequence of $\sigma$-fields
$$\{\Omega, \phi\} = \mathcal{F}_0 \subset \mathcal{F}_1 \subset \cdots \subset \mathcal{F}_n \subset \cdots \subset \mathcal{F}.$$
Let $H$ be a real or complex Hilbert space with norm $|\cdot|$. A sequence of $H$-valued strongly integrable functions $(f_n)_{n \geq 1}$ is a martingale if for each $n \geq 1$, $f_n$ is strongly measurable relative to $\mathcal{F}_n$, and for $n \geq 2$,
$$\E(d_n|\mathcal{F}_{n-1}) = 0 \quad \text{a.e.}$$
Here the difference sequence $(d_n)_{n \geq 1}$ is defined by $f_n = \sum_{i=1}^n d_i$, $n \geq 1$. In the following, we also call the limit $f = \sum\limits_{i=1}^\infty d_n$ martingale if the corresponding sequence $(f_n)_{n \geq 1}$ is a martingale.

A martingale is called conditionally symmetric if $d_{n+1}$ and $-d_{n+1}$ have the same conditional distribution given $d_1,...,d_n$.

We can write the Haar series as
 $$ g = (f, \varphi) \varphi + \sum_{k=0}^\infty \sum_{2^k\leq j \leq 2^{k+1}-1} (g,\Psi_{j,k})\Psi_{j,k},$$
 where
 \begin{eqnarray*}
  \varphi(x) &=& \begin{cases}
                  1,& x \in [0,1), \\
		  0,& \text{otherwise},
                 \end{cases} \\
  \Psi_{0,0}(x) &=& \varphi(2x) - \varphi(2x-1), \\
  \Psi_{j,k}(x) &=& 2^{k/2} \cdot \Psi_{0,0} (2^k x - j).
 \end{eqnarray*}
 Consider the probability space $\{\Omega, \mathcal{F}, \Prob\}$ defined by
 \begin{eqnarray*}
  \Omega &=& [0,1], \\
  \mathcal{F} &=& \mathcal{B}([0,1]), \\
  \Prob(A) &=& |A|, \quad A \in \mathcal{F},
 \end{eqnarray*}
and the sequence of $\sigma$-fields
\begin{eqnarray*}
 \{\Omega, \phi\} &=& \hspace{0.4cm} \{\Omega, \emptyset\} \\
		&& \subset \sigma(\varphi) \\
		&& \subset \sigma(\varphi, \Psi_{0,0}) \\
		&& \subset \sigma(\varphi, \Psi_{0,0}, \Psi_{2,1}) \\
		&& \subset \sigma(\varphi, \Psi_{0,0}, \Psi_{2,1},\Psi_{3,1}) \subset \cdots \\
		&& \subset \sigma(\varphi,\cdots,\Psi_{2^k,k},\cdots,\Psi_{2^{k+1}-1,k}, \Psi_{2^{k+1},k+1}) \subset \cdots \\
		&& \subset \mathcal{F}.
\end{eqnarray*}
We define
$$d_0 = (g,\varphi)\varphi, \quad d_1 = (g,\Psi_{0,0})\Psi_{0,0}, \quad d_2 = (g,\Psi_{2,1})\Psi_{2,1},\cdots$$
and
$$c_0 = (g,\varphi), \quad c_1 = (g,\Psi_{0,0}), \quad c_2 = (g,\Psi_{2,1}), \cdots $$
where the indices of $(f,\Psi_{\cdot,\cdot})\Psi_{\cdot,\cdot}$ and $(f,\Psi_{\cdot,\cdot})$ increase as in the definition of the sequence of $\sigma$-fields.

For each fixed $n \in \N_0$ and each $i=0,\cdots,n$, each of the sets $\{x:d_i(x)=c_i\}$, $\{x:d_i(x)=-c_i\}$, and $\{x:d_i(x)=0\}$ is either a superset of the support of $d_{n+1}$ or each of the sets and the support $d_{n+1}$ are disjoint. This implies that
\begin{eqnarray*}
 &&\Prob(d_{n+1} = c_{n+1} | d_i = j_i, \hspace{0.2cm} i \in \{1,...,n\}),\hspace{0.2cm} j_i \in \{c_i,-c_i,0\}) \\
 &=& \Prob(d_{n+1} = -c_{n+1} | d_i = j_i, \hspace{0.2cm} i \in \{1,...,n\}),\hspace{0.2cm} j_i \in \{c_i,-c_i,0\})
\end{eqnarray*}
so that the conditional distribution of $d_{n+1}$ and $-d_{n+1}$ is the same given $d_1,...,d_n$. Furthermore, we have $\E (d_{n+1}|\mathcal{F}_{n}) = 0$.
\end{proof}

\section{Extension of the calculation to the multidimensional case}

A very common way to extend the Haar basis to $[0,1]^d$ is given by the following construction (cf.~\cite{DeV98}). Let $E$ denote the collection of nonzero vertices of $[0,1]^d$. For each $e \in E$, we define the multivariate functions
$$\Psi^e(x_1,\cdots,x_d) := \Psi^{e_1}(x_1) \cdots \Psi^{e_d}(x_d),$$
where $\Psi^0(x) = \varphi(x)$, $\Psi^1(x) = \Psi_{0,0}(x)$. Furthermore, let 
$$\Psi_{j,k}^e(x) = 2^{kd/2}\cdot \Psi^e(2^k x - j), \quad k \geq 0, \quad 2^k \leq j_i \leq 2^{k+1}-1, \quad i=1,...,d$$
and
$$\Psi^*(x) = 1, \quad x \in [0,1]^d.$$
Then the collection of functions $\Psi^*$, $\Psi_{j,k}^e$, $e \in E$, $k \geq 0$, $2^k \leq j_i \leq 2^{k+1}-1$, $i=1,\cdots,d$ forms a basis for $L^p[0,1]^d$.

By considering the set $\mathcal{D}$ of dyadic cubes $I$ which form the supports of the functions $\Psi^*$, $\Psi_{j,k}^e$ and exchanging the notation of $\Psi^*$, $\Psi_{j,k}^e$ to $H_I$, we can also write the multivariate Haar basis as
$$\mathcal{H} = \{H_I\}_{I \in \mathcal{D}}.$$

\begin{lemma} \label{lemma:cond_mart2}
 Consider $f \in L^p[0,1]^d$ with corresponding Haar series
$$f =(f, \Psi^*) \Psi^* + \sum\limits_{e \in E} \sum\limits_{k=0}^\infty \sum\limits_{\substack{2^k \leq j_i \leq 2^{k+1}-1 \\ i=1,\cdots,d}} (f, \Psi_{j,k}^e) \Psi_{j,k}^e.$$

Then the inner double sum
$$ \sum\limits_{k=0}^\infty \sum\limits_{\substack{2^k \leq j_i \leq 2^{k+1}-1 \\ i=1,\cdots,d}} (f, \Psi_{j,k}^e) \Psi_{j,k}^e $$
forms a conditionally symmetric martingale on $[0,1]^d$ for each fixed $e \in E$, but so does not the Haar series itself.
\end{lemma}

\begin{proof}
 First we show that the Haar series itself does not form a conditionally symmetric martingale.

Let us assume that $d=2$ and remark that the proof goes analogously for $d > 2$. We have
\begin{eqnarray*}
  \Psi_{(0,0),0}^{(0,1)}(x_1,x_2) &=& \begin{cases}
                           1, & (x_1,x_2) \in [0,1] \times [0,\frac{1}{2}], \\
			   -1, & (x_1,x_2) \in [0,1] \times (\frac{1}{2},1], \\
			   0, &\text{otherwise},
                          \end{cases} \\
  \Psi_{(0,0),0}^{(1,0)}(x_1,x_2) &=& \begin{cases}
                           1, & (x_1,x_2) \in [0,\frac{1}{2}] \times [0,1],\\
			   -1, & (x_1,x_2) \in (\frac{1}{2},1] \times [0,1], \\
			   0, &\text{otherwise},
                          \end{cases} \\
  \Psi_{(0,0),0}^{(1,1)}(x_1,x_2) &=& \begin{cases}
                           1, & (x_1,x_2) \in [0,\frac{1}{2}] \times [0,\frac{1}{2}]\ \bigcup\ (\frac{1}{2},1] \times (\frac{1}{2},1], \\
			   -1, & (x_1,x_2) \in (\frac{1}{2},1] \times [0,\frac{1}{2}]\ \bigcup\ [0,\frac{1}{2}] \times (\frac{1}{2},1], \\
			   0, &\text{otherwise.}
                          \end{cases} \\
 \ 
\end{eqnarray*}
Therefore, the functions $\Psi_{(0,0),0}^{(0,1)}$, $\Psi_{(0,0),0}^{(1,0)}$, and $\Psi_{(0,0),0}^{(1,1)}$ can be represented as in Figure~\ref{fig:psi}.

\begin{figure}[htbp]
\includegraphics[width=14cm,height=6cm]{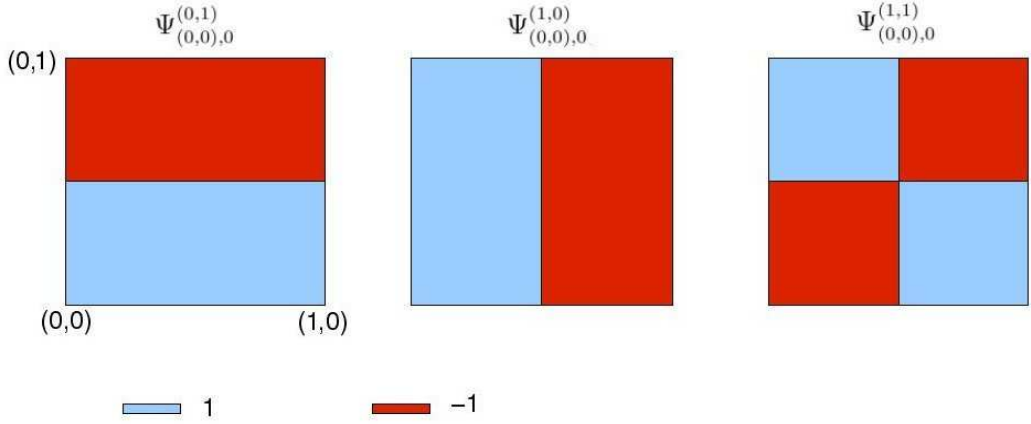}
\caption{The functions $\Psi_{(0,0),0}^{(0,1)}$, $\Psi_{(0,0),0}^{(1,0)}$, and $\Psi_{(0,0),0}^{(1,1)}$.}
\label{fig:psi}
\end{figure}

Thus
\begin{eqnarray*}
 \Prob(\Psi_{(0,0),0}^{(1,1)}=1|\Psi_{(0,0),0}^{(0,1)}= 1,\Psi_{(0,0),0}^{(1,0)}= 1) &=& 1, \\
 \Prob(\Psi_{(0,0),0}^{(0,1)}=1|\Psi_{(0,0),0}^{(1,0)}= 1,\Psi_{(0,0),0}^{(1,1)}= 1) &=& 1, \\
 \Prob(\Psi_{(0,0),0}^{(1,0)}=1|\Psi_{(0,0),0}^{(0,1)}= 1,\Psi_{(0,0),0}^{(1,1)}= 1) &=& 1, \\
 \end{eqnarray*}
 which implies that
\begin{eqnarray*}
 \E\left( \left. \Psi_{(0,0),0}^{(1,1)}\right|\Psi_{(0,0),0}^{(0,1)}= 1,\Psi_{(0,0),0}^{(1,0)}= 1\right) &=& 1, \\
 \E\left( \left. \Psi_{(0,0),0}^{(0,1)}\right|\Psi_{(0,0),0}^{(1,0)}= 1,\Psi_{(0,0),0}^{(1,1)}= 1\right) &=& 1, \\
 \E\left( \left. \Psi_{(0,0),0}^{(1,0)}\right|\Psi_{(0,0),0}^{(0,1)}= 1,\Psi_{(0,0),0}^{(1,1)}= 1\right) &=& 1. \\
\end{eqnarray*}

Therefore the multiparameter Haar series cannot a martingale.

Let us now consider the probability space $(\Omega, \mathcal{F}, \Prob)$ with
\begin{eqnarray*}
 \Omega &=& [0,1]^d, \\
 \mathcal{F} &=& \mathcal{B}([0,1]^d), \\
 \Prob(A) &=& |A|, \quad A \subset \mathcal{F},
\end{eqnarray*}
with the sequence of $\sigma$-fields
\begin{eqnarray*}
 \{\Omega, \phi\} &=& \hspace{0.4cm} \sigma(\Psi_{(0,0),0}^e) \\
		&&\subset \sigma(\Psi_{(0,0),0}^e, \Psi_{(2,2),1}^e) \subset \cdots \\
		&&\subset \sigma(\Psi_{(0,0),0}^e,\cdots, \Psi_{(3,2),1}^e, \Psi_{(2,3),1}^e, \Psi_{(3,3),1}^e) \subset \cdots \\
		&&\subset \sigma(\Psi_{(0,0),0}^e,\cdots, \Psi_{(3,3),1}^e, \Psi_{(4,4),2}^e, \Psi_{(5,4),2}^e, \Psi_{(6,4),2}^e, \Psi_{(7,4),2}^e) \subset \cdots \\
		&& \subset \sigma(\Psi_{(0,0),0}^e,\cdots, \Psi_{(7,4),2}^e, \Psi_{(4,5),2}^e, \Psi_{(5,5),2}^e, \Psi_{(6,5),2}^e) \subset \cdots \\
		&&\subset \mathcal{F}
\end{eqnarray*}
for each vertex $e \in E$. We denote this sequence of $\sigma$-fields by
$$\mathcal{F}_0 = \sigma(\Psi_{(0,0),0}^e), \quad \mathcal{F}_1 = \sigma(\Psi_{(0,0),0}^e, \Psi_{(2,2),1}^e), \quad \cdots. $$

Furthermore, for a fixed $e \in E$, we consider the partial sums of the inner sums of the Haar series
\begin{equation}
 \sum\limits_{k=0}^\infty \sum\limits_{\substack{2^k \leq j_i \leq 2^{k+1}-1 \\i=1,\cdots,d}} (f, \Psi_{j,k}^e) \Psi_{j,k}^e \label{eq:partial}
\end{equation}
and show that they form a conditionally symmetric martingale.

Let us denote the difference sequence $(d_n)_{n \geq 1}$ of $(\ref{eq:partial})$ by
$$d_0 = (f, \Psi_{(0,0),0}^e) \Psi_{(0,0),0}^e, \quad d_1 = (f, \Psi_{(2,2),1}^e)\Psi_{(2,2),1}^e, , \quad \cdots.$$

The corresponding coefficients are denoted by
$$c_0 = (f, \Psi_{(0,0),0}^e), \quad c_1 = (f, \Psi_{(2,2),1}^e), \quad \cdots.$$

For a fixed $k$ and $e \in E$, the support of the functions $\Psi_{j,k}^e$ is disjoint. Furthermore, for $l < k$ and $e \in E$, each of the sets $\{x:\Psi_{j,l}^e(x)=1\}$, $\{x:\Psi_{j,l}^e(x)=-1\}$, and $\{x:\Psi_{j,l}^e(x)=0\}$ is either a superset of the support of $\Psi_{j,k}^e$ or each of the sets and the support of $\Psi_{j,k}^e$ are disjoint.

This implies that $-d_{n+1}$ and $d_{n+1}$ have the same conditional distribution given $d_1,...,d_n$ and therefore
$$E(d_n|\mathcal{F}_{n-1}) = 0.$$
Thus, the partial sums of 
$$ \sum\limits_{k=0}^\infty \sum\limits_{\substack{2^k \leq j_i \leq 2^{k+1}-1 \\ i=1,\cdots,d}} (f, \Psi_{j,k}^e) \Psi_{j,k}^e$$
form a conditionally symmetric martingale.
\end{proof}

It is clear that the series $(f, \Psi^*) \Psi^* + \sum\limits_{k=0}^\infty \sum\limits_{\substack{2^k \leq j_i \leq 2^{k+1}-1 \\ i=1,\cdots,d}} (f, \Psi_{j,k}^e) \Psi_{j,k}^e$ is also a martingale.

\begin{theorem}
 Let $1 < p \leq 2$. Then for any $g \in L^p[0,1]^d$ we have
\begin{eqnarray*}
 &&\Vert g - G_m^{p} (g, \mathcal{H}) \Vert_p \\ &&\leq  \left(2+\frac{1}{\left(1-\left(\frac{1}{2}\right)^{d/p}\right)^2} \right) \left(\left(2^d-1\right) \left(\max\left(p,\frac{p}{p-1}\right)-1\right) \right)^2 \sigma_m(g)_p.
\end{eqnarray*}
\end{theorem}

\begin{proof}
Using Lemma \ref{lemma:cond_mart2}, we get an estimate for the upper bound in the Littlewood-Paley inequality by additionally applying the triangle inequality. Let $e^* \in E$. Then
\begin{eqnarray}
 &&\left\Vert (g, \Psi^*) \Psi^* + \sum_{e \in E} \sum_{k=0}^\infty \sum\limits_{\substack{2^k \leq j_i \leq 2^{k+1}-1 \\ i=1,\cdots,d}} (g, \Psi_{j,k}^e) \Psi_{j,k}^e \right\Vert_p \nonumber\\
 &\leq& \left\Vert (g, \Psi^*) \Psi^* + \sum_{k=0}^\infty \sum\limits_{\substack{2^k \leq j_i \leq 2^{k+1}-1 \\ i=1,\cdots,d}} (g, \Psi_{j,k}^{e^*}) \Psi_{j,k}^{e^*} \right\Vert_p \nonumber\\
 && + \sum_{e \in E\setminus \{e^*\}} \left\Vert\sum_{k=0}^\infty \sum\limits_{\substack{2^k \leq j_i \leq 2^{k+1}-1 \\ i=1,\cdots,d}} (g, \Psi_{j,k}^e) \Psi_{j,k}^e \right\Vert_p \nonumber\\
	&\leq& \sum\limits_{e \in E} C_4(p) \left\Vert \left( \sum_{I \in \mathcal{D}} |c_I(g) H_I|^2 \right)^{\frac{1}{2}} \right\Vert_p \nonumber\\
	&=& \left(2^d-1\right) C_4(p) \left\Vert \left( \sum_{I \in \mathcal{D}} |c_I(g) H_I|^2 \right)^{\frac{1}{2}} \right\Vert_p. \label{eq:littlewood2}
\end{eqnarray}

We now apply the method of duality (cf.~\cite{Lac08}) in order to determine the lower bound of the Littlewood-Paley inequality. The idea is to consider
$$S(g) := \left\Vert \left( \sum_{I \in \mathcal{D}} |c_I(g) H_I|^2 \right)^{\frac{1}{2}} \right\Vert_p $$
as an element of $L_{l^2(\mathcal{D})}^p$, that is
$$\varphi = \left\{\sqrt{|c_I(g) H_I|^2}: I \in \mathcal{D}\right\} $$
is considered as a $p$-integrable function taking values in $l^2(\mathcal{D})$. Due to the Hahn-Banach theorem, the dual function $\gamma = \{\gamma_I(x): I \in \mathcal{D}\} \in L_{l^2(\mathcal{D})}^{p'}$ is of norm one and satisfies
$$\Vert \varphi \Vert_{L_{l^2}^p} = (\varphi,\gamma) = \sum_{I \in \mathcal{D}} \sqrt{|c_I(g) H_I|^2} \int\limits_I \gamma_I dy, $$
where $p'$ is the conjugate index, i.~e. $1/p + 1/p' = 1$. This implies that we can assume that $\gamma_I$ is supported on $I$ and constant on $I$ since in the above formula, only the mean value of $\gamma_I$ over $I$ is important.

By defining the function
$$h := \sum\limits_{I \in \mathcal{D}} \left(\gamma_I \sqrt{|I|}\right) H_I,$$
we have on the one hand
$$S(h) = \Vert \gamma \Vert_{l^2(\mathcal{D})}$$
and on the other hand
\begin{eqnarray*}
 \Vert S(g) \Vert_p = \Vert \varphi \Vert_{L_{l^2}^p} &=& (\varphi,\gamma) \\
	&=& \sum\limits_{I \in \mathcal{D}} c_I(g) \gamma_I \sqrt{|I|} \\
	&=& (g,h) \\
	&\leq& \Vert g \Vert_p \Vert h \Vert_{p'} \\
	&\overset{\textup{(\ref{eq:littlewood2})}}{\leq}& \Vert g \Vert_p \left(2^d-1\right) C_4(p) \Vert S(h) \Vert_{p'} \\
	&=& \Vert g \Vert_p \left(2^d-1\right) C_4(p) \left\Vert \Vert \gamma \Vert_{l^2(\mathcal{D})} \right\Vert_{p'} \\
	&=& \Vert g \Vert_p \left(2^d-1\right) C_4(p).
\end{eqnarray*}
so that in the multidimensional case, the Littlewood-Paley inequality reads
\begin{equation}
 \frac{1}{C_4^*(p,d)} \left\Vert \left( \sum_I |c_I(g) H_I|^2 \right)^{\frac{1}{2}} \right\Vert_p \leq ||g||_p \leq C_4^*(p,d) \left\Vert \left( \sum_I |c_I(g) H_I|^2 \right)^{\frac{1}{2}} \right\Vert_p,
\end{equation}
where $C_4^*(p,d) = \left(2^d - 1\right) C_4(p) = \left(2^d - 1\right)\left(\max\left(p,\frac{p}{p-1}\right)-1\right)$.

Now, the remainder of the proof goes as for the univariate case. We note that
\begin{eqnarray*}
 \Vert g - S_\Lambda(g) \Vert_p &\leq& C_4^*(p,d)^2 \cdot \sigma_m(g)_p, \\
 \Vert S_{\Lambda_m \setminus \Lambda}(g) \Vert_p &\leq& C_4^*(p,d)^2 \cdot \sigma_m(g)_p,
\end{eqnarray*}
and
$$ \Vert S_{\Lambda \setminus \Lambda_m}(g) \Vert_p \leq \frac{1}{\left(1-\left(\frac{1}{2}\right)^{d/p}\right)^2} \Vert S_{\Lambda_m \setminus \Lambda}(g) \Vert_p$$
which will be derived in the following lemmas.

Combining the last three inequalities, we get
\begin{eqnarray*}
 	&&\Vert g - G_m^p(g)\Vert_p\\
 	&\leq& \Vert g - S_\Lambda(g) \Vert_p + \Vert S_\Lambda(g) - S_{\Lambda_m}(g) \Vert_p \\
	&=& \Vert g - S_\Lambda(g) \Vert_p + \Vert S_{\Lambda \setminus \Lambda_m}(g) - S_{\Lambda_m \setminus \Lambda}(g) \Vert_p \\
	&\leq& \Vert g - S_\Lambda(g) \Vert_p + \Vert S_{\Lambda \setminus \Lambda_m}(g)\Vert_p + \Vert S_{\Lambda_m \setminus \Lambda}(g) \Vert_p \\
	&\leq& \left(2+\frac{1}{\left(1-\left(\frac{1}{2}\right)^{d/p}\right)^2} \right) \left(\left(2^d-1\right) \left(\max\left(p,\frac{p}{p-1}\right)-1\right) \right)^2 \sigma_m(g)_p.
\end{eqnarray*}
\end{proof}

\begin{lemma} \label{lemma:5}
  Let $n_1 < n_2 < \cdots < n_s$ be integers and let $E_j \subset [0,1]^d$ be measurable sets, $j=1,...,s$. Then for any $0 < q < \infty$ we have
 \begin{equation}
  \int_{[0,1]^d} \left(\sum_{j=1}^s 2^{n_j d/q} \chi_{E_j}(x) \right)^q dx \leq \left(\frac{1}{1-\left(\frac{1}{2}\right)^{d/q}}\right)^q \cdot \sum_{j=1}^s 2^{n_j d} |E_j|. \nonumber
 \end{equation}
 \end{lemma}
 \begin{proof}
  Denote
	$$F(x) := \sum_{j=1}^s 2^{n_j d/q} \chi_{E_j}(x)$$
  and estimate it on the sets
  $$E_l^- := E_l \setminus \bigcup_{k=l+1}^s E_k, \quad l=1,...,s-1; \quad E_s^- := E_s.$$
  We have for $x \in E_l^-$
  \begin{eqnarray*}
   F(x) &\leq& \sum_{j=1}^l 2^{n_j d/q}
	= 2^{n_l d/q} \left( \frac{2^{n_1 d/q}}{2^{n_l d/q}} + \cdots +1\right) \\
	&\leq& 2^{n_l d/q} \sum_{i=0}^\infty \left(\frac{1}{2^{d/q}}\right)^i
	= 2^{n_l d/q} \frac{1}{1-\left(\frac{1}{2}\right)^{d/q}}.
  \end{eqnarray*}
  Therefore,
 $$ \int_0^1 F(x)^q dx \leq \left(\frac{1}{1-\left(\frac{1}{2}\right)^{d/q}} \right)^q \sum_{l=1}^s 2^{n_l d} |E_l^-| \leq \left(\frac{1}{1-\left(\frac{1}{2}\right)^{d/q}} \right)^q \sum_{l=1}^s 2^{n_l d} |E_l|,$$
which proves the lemma.
 \end{proof}
 
 \begin{lemma} \label{lemma:6}
  Consider
	$$f = \sum_{I \in Q} c_I H_I, \quad |Q| = N.$$
 Let $1 \leq p < \infty$. Assume that
 \begin{equation}
  \Vert c_I H_I \Vert_p \leq 1, \quad I \in Q. \label{eq:assumption2}
 \end{equation}

 Then
	$$\Vert f \Vert_p \leq \frac{1}{1-\left(\frac{1}{2}\right)^{d/p}} N^{1/p}.$$
 \end{lemma}

 \begin{proof}
  Denote by $n_1 < n_2 < \cdots < n_s$ all integers such that there is $I \in Q$ with $|I| = 2^{-d n_j}$. Introduce the sets
  $$E_j := \bigcup_{I \in Q: |I| = 2^{-d n_j}} I.$$
  Then the number $N$ of elements in $Q$ can be written in the form
  $$N = \sum_{j=1}^s |E_j| 2^{d n_j}.$$
  Furthermore, we have for $|I| = 2^{-kd}$, $k=0,1,2,\cdots$
  \begin{eqnarray*}
   \Vert c_i H_I\Vert_p = |c_I| \left(\int\limits_I \left|2^{dk/2}\right|^p dx\right)^{1/p} 
	=|c_I|\cdot 2^{dk/2} \cdot 2^{-kd/p} 
	= |c_I| |I|^{1/p-1/2}.
  \end{eqnarray*}
  The assumption (\ref{eq:assumption2}) implies
	$$ |c_I| \leq |I|^{1/2 - 1/p}.$$
  Next, we have
	$$\Vert f \Vert_p \leq \left\Vert \sum_{I \in Q} |c_I H_I| \right\Vert_p \leq \left\Vert \sum_{I \in Q} |I|^{-1/p} \chi_I(x) \right\Vert_p.$$
  The right hand side of this inequality can be rewritten as
  $$Y := \left( \int_{[0,1]^d} \left( \sum_{j=1}^s 2^{d n_j/p} \chi_{E_j}(x) \right)^p dx \right)^{1/p}.$$
  Applying Lemma \ref{lemma:5} with $q=p$, we get
  $$\Vert f \Vert_p \leq Y \leq \frac{1}{1-\left(\frac{1}{2}\right)^{d/p}} \left(\sum_{j=1}^s |E_j| 2^{d n_j} \right)^{1/p} = \frac{1}{1-\left(\frac{1}{2}\right)^{d/p}} N^{1/p}.$$
 \end{proof}

 \begin{lemma} \label{lemma:7}
  Consider
	$$f = \sum_{I \in Q} c_I H_I, \quad |Q| = N.$$
 Let $1 \leq p < \infty$. Assume that
	\begin{equation}
	 \Vert c_I H_I \Vert_p \geq 1, \quad I \in Q.
	\end{equation}
 Then
	$$\Vert f \Vert_p \geq \left(1-\left(\frac{1}{2}\right)^{d/p}\right) N^{1/p}.$$
 \end{lemma}
 \begin{proof}
  Define
  $$ u := \sum_{I \in Q} \overline{c}_I |c_I|^{-1} |I|^{1/p-1/2} H_I,$$
  where the bar means complex conjugate number. Then for $p' = \frac{p}{p-1}$ we have
  $$\Vert \overline{c}_I |c_I|^{-1} |I|^{1/p-1/2} H_I \Vert_{p'} = 1$$
  and, by Lemma \ref{lemma:5}
  $$\Vert u \Vert_{p'} \leq \frac{1}{1-\left(\frac{1}{2}\right)^{d/p}} N^{1/p'}.$$
  Consider $(f,u)$. We have on the one hand
  $$( f,u)  = \sum_{I \in Q} |c_I| |I|^{1/p-1/2} = \sum_{I \in Q} \Vert c_I H_I \Vert_p \geq N,$$
  and on the other hand
  $$ (f,u)  \leq \Vert f \Vert_p \Vert u \Vert_{p'},$$
  so that
  $$N \leq (f,u)  \leq \Vert f \Vert_p \Vert u \Vert_{p'} \leq \Vert f \Vert_p \frac{1}{1-\left(\frac{1}{2}\right)^{d/p}} N^{1/p'}$$
  which implies
  $$\Vert f \Vert_p \geq \left(1-\left(\frac{1}{2}\right)^{d/p}\right) N^{1/p}.$$
 \end{proof}

\begin{lemma}
 Let $1 < p < \infty$. Then for any $g \in L^p[0,1]^d$ we have
 $$\Vert S_{\Lambda \setminus \Lambda_m}(g) \Vert_p \leq \frac{1}{\left(1-\left(\frac{1}{2}\right)^{d/p}\right)^2} \cdot \Vert S_{\Lambda_m \setminus \Lambda}(g) \Vert_p.$$
\end{lemma}

\begin{proof}
 Denote
 $$A := \max\limits_{I \in \Lambda \setminus \Lambda_m} \Vert c_I(g) H_I \Vert_p \quad \text{and} \quad B:= \min\limits_{I \in \Lambda_m \setminus \Lambda} \Vert c_I(g) H_I \Vert_p.$$
 Then by the definition of $\Lambda_m$ we have
 $$ B \geq A.$$
 Using Lemma \ref{lemma:5}, we get
 \begin{equation}
  \Vert S_{\Lambda \setminus \Lambda_m}(g)\Vert_p \leq A \cdot \frac{1}{1-\left(\frac{1}{2}\right)^{d/p}} \cdot|\Lambda \setminus \Lambda_m |^{1/p} \leq B \cdot \frac{1}{1-\left(\frac{1}{2}\right)^{d/p}} \cdot|\Lambda \setminus \Lambda_m|^{1/p}. \label{eq:temp8}
 \end{equation}

 Using Lemma \ref{lemma:6}, we get
 $$\Vert S_{\Lambda_m \setminus \Lambda}(g)\Vert_p \geq B \cdot \left(1-\left(\frac{1}{2}\right)^{d/p}\right) \cdot |\Lambda_m \setminus \Lambda|^{1/p}$$
 so that
 \begin{equation}
  |\Lambda_m \setminus \Lambda|^{1/p} \leq \frac{1}{B \cdot \left(1-\left(\frac{1}{2}\right)^{d/p}\right)} \Vert S_{\Lambda_m \setminus \Lambda}(g)\Vert_p. \label{eq:temp9}
 \end{equation}
 Taking into account that $|\Lambda_m \setminus \Lambda| = |\Lambda \setminus \Lambda_m|$, we get
 \begin{eqnarray*}
  \Vert S_{\Lambda \setminus \Lambda_m}(g)\Vert_p &\overset{\textup{(\ref{eq:temp8})}}{\leq} B \cdot \frac{1}{1-\left(\frac{1}{2}\right)^{d/p}} \cdot|\Lambda \setminus \Lambda_m|^{1/p}
	\overset{\textup{(\ref{eq:temp9})}}{\leq}& \frac{1}{\left(1-\left(\frac{1}{2}\right)^{d/p}\right)^2} \Vert S_{\Lambda_m \setminus \Lambda}(g)\Vert_p.
 \end{eqnarray*}
\end{proof}

\end{document}